\documentclass[final,numberedheadings]{aipproc}

\layoutstyle{6x9}

\usepackage[all]{xy}
\usepackage{amsfonts}
\usepackage{graphicx}
\usepackage{amsmath}
\usepackage{amssymb}

\newtheorem{theorem}{Theorem}[section]

\begin{document}
\title[Some remarks on cosymplectic 3-structures]{Some remarks on cosymplectic 3-structures}



\author{Beniamino Cappelletti Montano}{
address={Universit\`{a} degli Studi di Cagliari, Dipartimento di Matematica e Informatica, Via Ospedale 72, 09124 Cagliari, Italia\\   b.cappellettimontano@gmail.com}}

\author{Antonio De Nicola}{
address={CMUC, Department of Mathematics, University of Coimbra, 3001-454 Coimbra, Portugal\\   antondenicola@gmail.com}}

\author{Ivan Yudin}{
address={CMUC, Department of Mathematics, University of Coimbra, 3001-454 Coimbra, Portugal\\   yudin@mat.uc.pt}}

\begin{abstract}
In this note we briefly review some recent results of the authors on the topological and geometrical properties of  $3$-cosymplectic manifolds.
\end{abstract}

\maketitle

\section{Introduction}
Cosymplectic manifolds were introduced in the frame of quasi-Sasakian manifolds by Blair in \cite{blairqs} as the closest odd-dimensional counterpart of
K\"{a}hler manifolds. Since then cosymplectic geometry has attracted the interest of many researchers also due to its role in mechanics and physics.
Recently, a great deal of work  on the topological properties of cosymplectic
manifolds  was done (see \cite{chineamarrero,chineamarrero2,li} among others).
In particular, in \cite{chineamarrero} Chinea, de Le\'{o}n and Marrero proved several important results for the Betti numbers of a
compact cosymplectic manifold.

The notion of $3$-cosymplectic manifold is the transposition of the notion of
cosymplectic manifold to the setting of $3$-structures. Namely, a $3$-cosymplectic
manifold is a smooth manifold endowed with three distinct cosymplectic
structures related to each other by means of some relations formally similar to the
quaternionic identities (see Section~\ref{preliminaries} for more details). This note contains a concise review of the main properties of $3$-cosymplectic
manifolds, recently obtained by the authors in \cite{cappellettidenicola,
cappellettidenicolayudin}. Especially, we emphasize our results
concerning Betti numbers of compact $3$-cosymplectic manifolds. Finally, we
present a method for constructing non-trivial examples of such compact
manifolds.

\section{3-cosymplectic geometry}\label{preliminaries}
An \emph{almost contact manifold} is an odd-dimensional smooth manifold $M$
endowed with a tensor field $\phi$ of endomorphisms on the tangent spaces, a
vector field $\xi$ and a $1$-form $\eta$ satisfying \ $\phi^{2}=-I+\eta\otimes\xi$, \ where $I$ denotes the identity mapping of $TM$. It is known that
there exists a Riemannian metric $g$ which is compatible with the structure, in the sense that
\begin{equation}\label{metric1}
g(\phi X, \phi Y) = g(X,Y) - \eta(X)\eta(Y)
\end{equation}
for any $X,Y\in\Gamma(TM)$. When one fixes one compatible metric, the resulting
geometric structure $(\phi,\xi,\eta,g)$ is called an \emph{almost
contact metric structure} on $M$. From \eqref{metric1} it follows that the
bilinear form  $\Phi:=g(\cdot,\phi\cdot)$ is in fact a $2$-form,
called the \emph{fundamental $2$-form} of the almost contact metric manifold. An \emph{almost cosymplectic manifold} is an almost contact metric manifold
$(M,\phi,\xi,\eta,g)$ such that both the $1$-form $\eta$ and the fundamental $2$-form $\Phi$ are closed. If in addition the structure is \emph{normal},
that is, if the Nijenhuis  tensor field of $\phi$ vanishes identically, $(M,\phi,\xi,\eta,g)$ is said to be a \emph{cosymplectic manifold}. In terms of the
covariant derivative of the structure tensor field $\phi$, this condition is
equivalent to  $\nabla\phi=0$.   Now, we come to the main topic of the paper.
A triple of almost contact structures $(\phi_1,\xi_1,\eta_1)$, $(\phi_2,\xi_2,\eta_2)$,
$(\phi_3,\xi_3,\eta_3)$ on a manifold $M$, related by the identities
\begin{equation} \label{quaternionic}
\begin{split}
\phi_\gamma=\phi_{\alpha}\phi_{\beta}-\eta_{\beta}\otimes\xi_{\alpha}=-\phi_{\beta}\phi_{\alpha}+\eta_{\alpha}\otimes\xi_{\beta},\quad\\
\xi_{\gamma}=\phi_{\alpha}\xi_{\beta}=-\phi_{\beta}\xi_{\alpha}, \ \ \eta_{\gamma}= \eta_{\alpha}\circ\phi_{\beta}=-\eta_{\beta}\circ\phi_{\alpha},
\end{split}
\end{equation}
for any even permutation $(\alpha,\beta,\gamma)$ of the set
$\left\{1,2,3\right\}$, is called an \emph{almost contact 3-structure} on
$M$.
Then, the dimension of the manifold is necessarily of the form $4n+3$. This notion was introduced independently by Kuo (\cite{kuo}) and Udriste
(\cite{udriste}). In particular, Kuo proved that one can always find a Riemannian metric $g$ which is compatible with each almost contact structure. If we
fix a compatible metric, we speak of \emph{almost contact metric 3-structure}. Any smooth manifold endowed with an almost contact metric $3$-structure carries two
orthogonal distributions: the \emph{Reeb distribution} ${\mathcal
V}:=\operatorname{span}\{\xi_1,\xi_2,\xi_3\}$ and the \emph{horizontal distribution} ${\mathcal
H}:=\ker(\eta_1)\cap\ker(\eta_2)\cap\ker(\eta_3)$.

A remarkable case is when each  structure is cosymplectic. In this case we say that $M$ is a \emph{3-cosymplectic manifold}. In any
$3$-cosymplectic manifold the forms $\eta_\alpha$ and $\Phi_\alpha$ are harmonic. Moreover, the tensors $\xi_\alpha$, $\eta_\alpha$, $\phi_\alpha$,
$\Phi_\alpha$ are all $\nabla$-parallel. In particular, since the Reeb vector fields commute with each other, it follows that the Reeb distribution is
integrable and defines a $3$-dimensional foliation ${\mathcal F}_{3}$ of $M$. As it was proven in \cite{cappellettidenicola}, ${\mathcal F}_{3}$ is a
Riemannian and transversely hyper-K\"{a}hler foliation with totally geodesic leaves. Moreover, since $d\eta_\alpha=0$, also the horizontal distribution
$\mathcal H$ is integrable and hence defines a Riemannian, totally geodesic foliation complementary to ${\mathcal F}_3$.

Another important property of $3$-cosymplectic manifolds that should be mentioned is that they are Ricci-flat (\cite{cappellettidenicola}).

\section{The cohomology of a $3$-cosymplectic manifold}

Let $M$ be a compact $3$-cosymplectic manifold of dimension $4n+3$. We will denote by $H^{\ast}_{dR}(M)$ the usual de Rham cohomology of $M$. By the
Hodge-de Rham theory each vector space $H^{k}_{dR}(M)$ can be identified with the vector space $\Omega^{k}_{H}(M)$ of harmonic $k$-forms on $M$.  Recall
also that the  space of basic $k$-forms  (with respect to ${\mathcal F}_3$) is defined by
\begin{equation*}
\Omega^{k}_{B}(M):=\left\{ \omega\in\Omega^{k}(M)\  \middle|\  i_{\xi_\alpha}\omega=0, \
i_{\xi_\alpha}d\omega=0, \ \hbox{ for each } \alpha=1,2,3 \right\}.
\end{equation*}
Since the differential $d$ preserves basic forms, it induces a cohomology $H^{\ast}_{B}(M)$ which is called \emph{basic cohomology}.

For each $\alpha\in\left\{1,2,3\right\}$ we define two linear operators $l_{\alpha}:\Omega^{k}(M)\longrightarrow\Omega^{k+1}(M)$,
$\omega\mapsto\eta_{\alpha}\wedge\omega$, and $\lambda_{\alpha}:\Omega^{k+1}(M)\longrightarrow\Omega^{k}(M)$, $\omega\mapsto i_{\xi_\alpha}\omega$.
Moreover, we define $e_\alpha:=l_{\alpha}\circ\lambda_{\alpha}$. By \cite[Proposition 1]{chineamarrero} the operators $l_\alpha$, $\lambda_\alpha$, and
hence $e_\alpha$, preserve harmonic forms. Then one can prove the following decomposition
\begin{equation}\label{decomposition}
\Omega^{k}_{H}(M)=\bigoplus_{\epsilon_{1}, \epsilon_{2}, \epsilon_{3}\in\{0,1\}}\Omega^{k}_{H,\epsilon_{1}\epsilon_{2}\epsilon_{3}}(M),
\end{equation}
where we have put, for each triple $\epsilon_{1},\epsilon_{2},\epsilon_{3}\in\left\{0,1\right\}$,
\begin{equation*}
\Omega^{k}_{H,\epsilon_{1}\epsilon_{2}\epsilon_{3}}(M):=\left\{\omega\in\Omega^{k}_{H}(M)\
\middle|\  e_{\alpha}\omega=\epsilon_{\alpha}\omega, \ \alpha=1,2,3\right\}.
\end{equation*}
Moreover, one can prove that the operators $l_1$, $l_2$, $l_3$ induce isomorphisms between the vector spaces
$\Omega^{\ast}_{H,\epsilon_{1}\epsilon_{2}\epsilon_{3}}$ according to the following diagram
\begin{equation*}
\xymatrix{ & \Omega^{k+1}_{H,100}\left( M \right) \ar[rr]^{l_2}
\ar[dd]^(.3){l_3}|!{[dl];[dr]}\hole &&
\Omega^{k+2}_{H,110} \left( M \right) \ar[dd]_{l_3}\\
\Omega^k_{H,000}\left( M \right)
\ar[rr]^(.7){l_2}\ar[dd]_{l_3}\ar[ru]^{l_1} &&
\Omega^{k+1}_{H,010}\left( M
\right)\ar[dd]^(.7){l_3} \ar[ru]^{l_1}\\
& \Omega^{k+2}_{H,101}\left( M \right)
{\ar[rr]^(.4){l_2}|!{[dr];[ur]}\hole}   &&
\Omega^{k+3}_{H,111}\left( M \right)\\
\Omega^{k+1}_{H,001}\left( M \right) \ar[rr]^{l_2}\ar[ru]^{l_1} &&
\Omega^{k+2}_{H,011}\left( M \right) \ar[ru]_{l_1} }
\end{equation*}
for each $0\le k\le 4n$. Therefore, the whole information about cohomology groups of $M$ is contained in the vector spaces $\Omega^k_{H,000}(M)$, $0\le k\le
4n$. It is worth to mention that $\Omega^k_{H,000}(M)$  can be identified with
the space of basic harmonic $k$-forms on $M$ (with respect to $\mathcal{F}_3$). In particular,
$b^{h}_{k}:=\dim(\Omega^k_{H,000}(M))$ is the $k$-th basic Betti number. Now, taking the decomposition \eqref{decomposition} into account and using the
above isomorphisms between the vector spaces $\Omega^{\ast}_{H,\epsilon_{1}\epsilon_{2}\epsilon_{3}}(M)$, one gets the following formula for the $k$-th
Betti number of $M$
\begin{equation}\label{betti1}
b_{k}=b^{h}_{k}+3b_{k-1}^{h}+3b_{k-2}^{h}+b_{k-3}^{h}.
\end{equation}
 On the other hand, one can prove (see \cite{cappellettidenicolayudin} for more details) that, for each odd integer $k$,
$\Omega^{k}_{H,000}(M)$ is a $\mathbb{H}$-module and thus $b^{h}_{k}$ is divisible by $4$. Then by \eqref{betti1} it follows that, for any odd integer $k$,
$b_{k-1}+b_{k}$ is divisible by $4$. Another restriction on the Betti numbers of a compact 3-cosymplectic manifold is the following inequality
\begin{equation}\label{betti2}
b_{k}\geq\binom{k+2}{2}
\end{equation}
for $0\leq k\leq 2n+1$, which is stronger than the analogous inequality for hyper-K\"{a}hler manifolds, due to Wakakuwa (\cite{wakakuwa}), namely
$b_{2k}\geq\binom{k+2}{2}$.  \ We conclude the section by describing
 an action of the Lie algebra $so(4,1)$ on $\Omega^{k}_{H,000}(M)$. For every even permutation $(\alpha,\beta,\gamma)$ of $\{1,2,3\}$ let us consider
the $2$-form $\Xi_\alpha:=\frac{1}{2}(\Phi_{\alpha}+2\eta_{\beta}\wedge\eta_{\gamma})$. Then we define the operators
$L_{\alpha}:\Omega^{k}(M)\longrightarrow\Omega^{k+2}(M)$ and $\Lambda_{\alpha}:\Omega^{k+2}(M)\longrightarrow\Omega^{k}(M)$ by
$L_{\alpha}\omega:=\Xi_{\alpha}\wedge\omega$ and $\Lambda_{\alpha}:=\ast L_{\alpha} \ast$. Since $L_\alpha$ and $\Lambda_\alpha$ preserve harmonicity, one
can consider them as endomorphisms  of $\Omega^{\ast}_{H,000}(M)$. Then, by \cite[Proposition 4.3]{cappellettidenicolayudin} one has that, on
$\Omega^{\ast}_{H,000}$, $[L_{\alpha},\Lambda_{\alpha}]=-H$, where $H:\Omega^{k}_{H,000}(M)\longrightarrow\Omega^{k}_{H,000}(M)$ is the operator defined by
$H\omega=(2n-k)\omega$. Moreover, for each $\alpha\in\{1,2,3\}$ we define another operator $K_\alpha$ on $\Omega^{k}_{H,000}(M)$ by
$K_{\alpha}:=[L_\beta,\Lambda_\gamma]$, where $(\alpha,\beta,\gamma)$ is an even permutation of $\{1,2,3\}$. Then we have the following result.

\begin{theorem}[\cite{cappellettidenicolayudin}]
The linear span $\mathfrak{g}$ of the operators $H$, $L_\alpha$, $\Lambda_\alpha$, $K_\alpha$, $\alpha\in\left\{1,2,3\right\}$, is a Lie algebra isomorphic to $so(4,1)$. Consequently $\Omega^{\ast}_{H,000}(M)$ is an $so(4,1)$-module.
\end{theorem}

\section{Examples of compact 3-cosymplectic manifolds}
The standard example of $3$-cosymplectic manifold is $\mathbb{R}^{4n+3}$ with the almost contact metric $3$-structure described in
\cite{cappellettidenicola} in terms of Darboux coordinates. Since this structure is invariant by translations, we get a $3$-cosymplectic structure on the
flat torus $\mathbb{T}^{4n+3}$ (see also \cite{martincabrera}). Both these
examples are global products of a hyper-K\"{a}hler manifold with an abelian
Lie group. In fact, locally this is always true: every $3$-cosymplectic manifold is locally a Riemannian product of a hyper-K\"{a}hler factor with a
$3$-dimensional flat abelian Lie group. Thus it makes sense to ask whether there
are
examples of $3$-cosymplectic manifolds which are not global products of
a hyper-K\"{a}hler manifold with a $3$-dimensional Lie group. The answer to this
question is affirmative and now we describe a procedure for constructing
such examples. \ Let $(M^{4n},J_{\alpha},G)$ be a compact hyper-K\"{a}hler manifold and $f$ a hyper-K\"{a}hler isometry on it. We define an action
$\varphi$ of $\mathbb{Z}^3$ on $M^{4n}\times\mathbb{R}^{3}$ by
\begin{equation*}
\varphi((k_1,k_2,k_3),(x,t_1,t_2,t_3)) = (f^{k_1+k_2+k_3}(x), t_{1}+k_{1},t_{2}+k_{2},t_{3}+k_{3}).
\end{equation*}
We define  a $3$-cosymplectic structure on the orbit space $M^{4n+3}_{f}:=(M^{4n}\times\mathbb{R}^{3})/\mathbb{Z}^{3}$ in the following way. Let us
consider the vector fields $\xi_\alpha:=\frac{\partial}{\partial t_\alpha}$ and the $1$-forms $\eta_\alpha:=d t_\alpha$ on $M^{4n}\times \mathbb{R}^{3}$.
Next we define, for each $\alpha\in\left\{1,2,3\right\}$, a tensor field $\phi_\alpha$ on $M^{4n}\times\mathbb{R}^{3}$ by putting
$\phi_{\alpha}X:=J_{\alpha}X$ for any $X\in\Gamma(T{M^{4n}})$ and $\phi_{\alpha}\xi_{\alpha}:=0$,
$\phi_{\alpha}\xi_{\beta}:=\epsilon_{\alpha\beta\gamma}\xi_{\gamma}$, where $\epsilon_{\alpha\beta\gamma}$ denotes the sign of the permutation
$(\alpha,\beta,\gamma)$ of $\left\{1,2,3\right\}$. Then
$(\phi_\alpha,\xi_\alpha,\eta_\alpha,g)$, where  $g$ denotes the product metric,
is a
$3$-cosymplectic structure on $M^{4n}\times\mathbb{R}^{3}$. Being invariant
under the action $\varphi$, the structure $(\phi_\alpha,\xi_\alpha,\eta_\alpha,g)$ descends
to a $3$-cosymplectic structure on $M^{4n+3}_{f}$. By using this general procedure we can construct non-trivial
examples of compact $3$-cosymplectic manifolds. In fact, let us consider the hyper-K\"{a}hler manifold $\mathbb{T}^{4}=\mathbb{H}/\mathbb{Z}^4$ and the
hyper-K\"{a}hler isometry $f$ given by the multiplication by the quaternionic unit
$\bold{i}$ on the right. Then $M^{7}_{f}:=(\mathbb{T}^{4}\times\mathbb{R}^{3})/\mathbb{Z}^{3}$, endowed
with the geometric structure described above, is a compact $3$-cosymplectic manifold which is not the global product of a compact $4$-dimensional
hyper-K\"{a}hler manifold $K^{4}$ with the flat torus. Indeed, we have only two possibilities for a compact $4$-dimensional hyper-K\"{a}hler manifold:
either $K^{4}\cong\mathbb{T}^{4}$ or it is a complex K3-surface. In the first case $b_{2}(K^{4}\times\mathbb{T}^{3})=21$, in the second
$b_{2}(K^{4}\times\mathbb{T}^{3})=25$. However, in \cite{cappellettidenicolayudin} it was proven that $b_{2}(M^{7}_{f})<21$.

Other examples can be obtained from the previous ones by applying a \emph{${\mathcal D}_{a}$-homothetic deformation}, that is a change of the structure
tensors of the following type
\begin{equation}\label{deformation}
\bar\phi:=\phi, \ \ \ \bar\xi:=\frac{1}{a}\xi, \ \ \ \bar\eta:=a \eta, \ \ \ \bar{g}:=a g + a(a-1)\eta\otimes\eta,
\end{equation}
where $a>0$. This notion was introduced by Tanno (\cite{tanno68}) in the contact
metric case, but it can be easily extended to the more general context of
almost contact metric structures. In particular, it can be proved that the class of cosymplectic structures is preserved by $\mathcal D$-homothetic
deformations. Now, let $(M,\phi_{\alpha},\xi_{\alpha},\eta_{\alpha},g)$, $\alpha\in\left\{1,2,3\right\}$, be a 3-cosymplectic manifold. Then by  applying
the same ${\mathcal D}_{a}$-homothetic deformation to each cosymplectic structure $(\phi_{\alpha},\xi_{\alpha},\eta_{\alpha},g)$, one obtains three new
cosymplectic structures $(\bar{\phi}_{1},\bar{\xi}_{1},\bar{\eta}_{1},\bar{g})$, $(\bar{\phi}_{2},\bar{\xi}_{2},\bar{\eta}_{2},\bar{g})$,
$(\bar{\phi}_{3},\bar{\xi}_{3},\bar{\eta}_{3},\bar{g})$, which are still related
to each other by means of the quaternionic-like relations
\eqref{quaternionic}. Thus $(\bar\phi_{\alpha},\bar\xi_{\alpha},\bar\eta_{\alpha},\bar{g})$, $\alpha\in\left\{1,2,3\right\}$, is a new 3-cosymplectic
structure on  $M$.  In particular, this procedure allows to define other $3$-cosymplectic structures on $M^{4n+3}_{f}$ from the structure described before.

We conclude with the following remark concerning the existence of $3$-cosymplectic structures on almost cosymplectic Einstein manifolds. In the context of
Sasakian manifolds, Apostolov, Draghici and Moroianu proved the following theorem:

\begin{theorem}[\cite{apostolov06}]\label{goldberg}
Let $(M,\phi,\xi,\eta,g)$ be a Sasakian Einstein manifold. Then any contact metric structure $(\phi',\xi',\eta',g)$ on $M$, with the same metric $g$ is
Sasakian. Moreover, if $\xi'\neq\pm\xi$, then either $(M,g)$ admits a 3-Sasakian structure or $(M,g)$ is covered by a round sphere.
\end{theorem}

It could be interesting to investigate on the cosymplectic counterpart, if any,
of Theorem \ref{goldberg}. In fact, this would permit to construct new
examples of 3-cosymplectic manifolds. In this context we mention the following result on compact Einstein cosymplectic manifolds.

\begin{theorem}[\cite{pastore10}]
Every  compact   Einstein almost  cosymplectic manifold $(M,\phi,\xi,\eta,g)$,
such that $\xi$ is Killing, is cosymplectic and Ricci-flat. Furthermore, any
other almost cosymplectic structure $(\phi',\xi',\eta',g)$  on $M$ is necessarily cosymplectic and Ricci-flat.
\end{theorem}

Notice that the proof of Theorem \ref{goldberg} does not work in the case of
cosymplectic manifolds. In fact, it  uses a property of the cone metric which
holds only for Sasakian manifolds.

\begin{theacknowledgments}
The second author is partially supported by the FCT grant PTDC/MAT/099880/2008 and by MICIN (Spain) grant MTM2009-13383. The last author is supported by
the FCT Grant SFRH/BPD/31788/2006.
\end{theacknowledgments}

\end{document}